\documentclass{amsart}
\usepackage{latexsym}
\usepackage{enumerate}
\usepackage{amsmath, amssymb, amsfonts, amsxtra, amsthm}
\usepackage[active]{srcltx}

\newcommand {\rea}{\mathbb{R}}
\newcommand {\BBZ}{\mathbb{Z}}
\newcommand {\calM}{\mathcal{M}}
\newcommand {\calS}{\mathcal{S}}
\newcommand {\calI}{\mathcal{I}}

\newtheorem{theorem}{Theorem}[section]

\newtheorem{lemma}[theorem]{Lemma}

\author{Richard Oberlin}
\address{
Mathematics Department\\
Louisiana State University\\
Baton Rouge, LA}
\email{oberlin@math.lsu.edu}
\thanks{The author is supported in part by NSF Grant DMS-1068523.}


\begin{document}
\title{A Marcinkiewicz maximal-multiplier theorem}

\begin{abstract}
For $r < 2$, we prove the boundedness of a maximal operator formed by applying all multipliers $m$ with $\|m\|_{V^r} \leq 1$ to a given function. 
\end{abstract}

\maketitle

\section{Introduction}
Given an exponent $r$ and a function $f$ defined on $\rea$, consider the $r$-variation norm
\[
\|f\|_{V^r} = \|f\|_{L^\infty} + \sup_{N, \xi_0 < \cdots < \xi_N} \left(\sum_{i = 1}^N|f(\xi_{i}) - f(\xi_{i - 1})|^r \right)^{1/r}
\]
where the supremum is over all strictly increasing finite length sequences of real numbers.

The classical Marcinkiewicz multiplier theorem states that if $r=1$ and a function $m$ is of bounded $r$-variation uniformly on dyadic shells, then $m$ is an $L^p$ multiplier for $1 < p < \infty$ and
\begin{equation} \label{mbound}
\|(m\hat{f})\check{\ }\|_{L^p} \leq C_{p,r} \sup_{k \in \BBZ} \|1_{D_k} m\|_{V^r}\|f\|_{L^p}
\end{equation}
where $D_k = [-2^{k+1}, -2^k) \cup (2^k,2^{k+1}]$ and $\hat{\ },\check{\ }$ denote the Fourier-transform and its inverse. Later, Coifman, Rubio de Francia, and Semmes \cite{coifman88mdf}, see also \cite{tao01emt}, showed that the requirement of bounded $1$-variation can be relaxed to allow for functions of bounded $2$-variation, and in fact \eqref{mbound} holds whenever $r \geq 2$ and $|\frac{1}{p} - \frac{1}{2}| < \frac{1}{r}.$

The estimate \cite{coifman88mdf} does not discriminate between multipliers of bounded $2$-variation and those of bounded $r$-variation where $r < 2$, and so one might ask whether there is anything to be gained by controlling the variation norm of multipliers in the latter range of exponents. 

Defining the maximal-multiplier operator
\begin{equation} \label{maxdef}
\calM_r[f](x) = \sup_{m : \|m\|_{V^r} \leq 1}|(m\hat{f})\check{\ }(x)| 
\end{equation}
where the supremum is over all functions in the $V^r$ unit ball, we have
\begin{theorem} \label{mainthm}
Suppose $1 \leq r < 2$ and $r < p < \infty.$ Then
\begin{equation} \label{mmbound}
\|\calM_r[f]\|_{L^p} \leq C_{p,r} \|f\|_{L^p}.
\end{equation}
\end{theorem}  
\noindent The case $r=1$ was observed independently by Lacey \cite{laceypc}.

Note that, in the definition of $\calM_r$, each $m$ is required to have finite $r$-variation on all of $\rea$ rather than simply on each dyadic shell as in \eqref{mbound}. This is necessary for boundedness, as can be seen from the counterexamples of Christ, Grafakos, Honz\'{i}k and Seeger \cite{christ05mfa}.

Although the maximal operator \eqref{maxdef} would seem to be fairly strong, we do not yet know of an application for the bound above. We will, however, quickly illustrate a strategy for its use that falls an (important) $\epsilon$ short of success. Let $\Psi$ be (say) a Schwartz function, and for each $\xi,x \in \rea$ and $k \in \BBZ$ consider the $2^k$-truncated partial Fourier integral
\[
\calS_k[f](\xi,x) = p.v. \int f(x - t) e^{2 \pi i \xi t} \Psi(2^{-k}t)\frac{1}{t}\ dt.
\]
It was proven by Demeter, Lacey, Tao, and Thiele \cite{demeter08bdr} that for $q = 2$ and $1 < p < \infty$
\begin{equation} \label{mmc}
\|\sup_{\|g\|_{L^q}=1}\|\sup_k|(\calS_k[f](\cdot,x)\hat{g})\check{\ }|\|_{L^q}\|_{L^p_x} \leq C_{p,q} \|f\|_{L^p}.
\end{equation}
If we had the bound
\begin{equation} \label{vtc}
\|\calS_k[f](\xi,x)\|_{L^p_x(\ell^{\infty}_k(V^r_\xi))} \leq C_{p,r} \|f\|_{L^p}
\end{equation}
for some $r < 2$, then an application of Theorem \ref{mainthm} would give \eqref{mmc} for $q > r$ by rather different means than \cite{demeter08bdr}. In fact, one can see by applying the method in Appendix D of \cite{oberlin09vnc} that \eqref{vtc} holds for $r > 2$ and $p > r'.$ Unfortunately, it does fail for $r \leq 2.$

\section{Proof of Theorem \ref{mainthm}}

The following lemma was proven in \cite{coifman88mdf}, see also \cite{lacey07irt}. 
\begin{lemma} \label{intervaldecompositionlemma}
Let $m$ be a compactly supported function on $\rea$ of bounded $r$-variation for some $1 \leq r < \infty.$ Then for each integer $j \geq 0$, one can find a collection $\Upsilon_j$ of pairwise disjoint subintervals of $\rea$ and coefficients $\{b_{\upsilon}\}_{\upsilon \in \Upsilon_j} \subset \rea$ so that $|\Upsilon_j| \leq 2^j$, $|b_{\upsilon}| \leq 2^{-j/r}\|m\|_{V_r} $, and
\begin{equation} \label{intervaldecompeq}
m = \sum_{j \geq 0} \sum_{\upsilon \in \Upsilon_j} b_{\upsilon} 1_{\upsilon} 
\end{equation}
where the sum in $j$ converges uniformly.
\end{lemma} 
The lemma above was applied in concert with Rubio de Francia's square function estimate \cite{rubiodefrancia85lpi} to obtain \eqref{mbound}. Here, we will argue similarly, exploiting the analogy between the Rubio de Francia square function estimate and the variation-norm Carleson theorem. 

It was proven in \cite{rubiodefrancia85lpi} that for $p \geq 2$
\[
\sup_{\calI}\|\big(\sum_{I \in \calI} |(1_I \hat{f})\check{\ }|^2 \big)^{1/2}\|_{L^p} \leq C_p \|f\|_{L^p}
\]
where the supremum above is over all collections of pairwise disjoint subintervals of $\rea$. Consider the partial Fourier integral
\[
\calS[f](\xi,x) = (1_{(\infty,\xi]}\hat{f})\check{\ }(x).
\]
It was proven in \cite{oberlin09vnc} that for $s > 2$ and $p > s'$
\[
\|\calS[f](\xi,x)\|_{L^p_x(V^s_\xi)} \leq C_{p,s} \|f\|_{L^p} 
\]
or, equivalently,
\begin{equation} \label{vncrub}
\|\sup_{\calI} \big(\sum_{I \in \calI} |(1_I \hat{f})\check{\ }|^s \big)^{1/s}\|_{L^p} \leq C_{p,s} \|f\|_{L^p}.
\end{equation}

Note that, by standarding limiting arguments, taking the supremum in \eqref{maxdef} to be over all compactly supported $m$ such that $\|m\|_{V^r} \leq 1$ does not change the definition of $\calM_r.$ Applying the decomposition \eqref{intervaldecompeq}, we see that for any compactly supported $m$ with $\|m\|_{V^r} \leq 1$ we have
\begin{align*}
|(m \hat{f})\check{\ }(x)|  &\leq \sum_{j \geq 0} \sum_{\upsilon \in \Upsilon_j} |b_{\upsilon} (1_{\upsilon}\hat{f})\check{\ }(x)| \\
&\leq \sum_{j \geq 0}\sup_{\upsilon \in \Upsilon_j} |b_\upsilon| |\Upsilon_j|^{\frac{1}{s'}} \big(\sum_{\upsilon \in \Upsilon_j}|(1_\upsilon \hat{f})\check{\ }(x)|^s \big)^{1/s} \\
& \leq C_{r,s} \sup_{\calI} \big(\sum_{I \in \calI} |(1_I \hat{f})\check{\ }(x)|^s \big)^{1/s}
\end{align*}
where, for the last inequality, we require $s < r'.$ 

Provided that $r < 2$ and $p > r$ we can choose an $s < r'$ with $s > 2$ and $p > s'$, giving \eqref{vncrub} and hence \eqref{mmbound}.   

The argument of Lacey \cite{laceypc} for $r=1$ follows a similar pattern, except with Marcinkiewicz's method in place of \cite{coifman88mdf} and the standard Carleson-Hunt theorem in place of \cite{oberlin09vnc}.

\bibliographystyle{hamsplain}
\bibliography{roberlin}
\end{document}